\documentclass[leqno,12pt,a4paper]{article}
\usepackage[T1]{fontenc}
\usepackage{fourier}
\usepackage{setspace} \usepackage{etoolbox} \usepackage{amsfonts,amssymb}
\usepackage[margin=1in]{geometry}
\usepackage{microtype}
\usepackage{times}
\usepackage[amsthm,thmmarks]{ntheorem}
\usepackage{amsmath,xypic,thmtools}
\usepackage{multicol}
\usepackage{hyphenat}
\usepackage[usenames,dvipsnames]{color}
\usepackage{rotating}
\usepackage{longtable}
\usepackage{pdflscape}
\usepackage{caption}
\usepackage{subcaption} \usepackage{tcolorbox}
\usepackage{stackengine}

\usepackage{epigraph}
\usepackage{enumitem}
\setlist[enumerate]{listparindent=0.5in}
\usepackage[cal=boondoxo,calscaled=.96,scr=rsfs]{mathalpha}
\DeclareMathAlphabet{\mathscrbf}{OMS}{mdugm}{b}{n} \newcommand{\be}{\begin{equation}}
\newcommand{\ee}{\end{equation}}
\newcommand{\bes}{\begin{equation*}}
\newcommand{\ees}{\end{equation*}}
\newcommand{\bea}{\begin{eqnarray}}
\newcommand{\eea}{\end{eqnarray}}
\newcommand{\beas}{\begin{eqnarray}}
\newcommand{\eeas}{\end{eqnarray}}
\newcommand{\ben}{\begin{note}}
\newcommand{\een}{\end{note}}
\newcommand{\bexl}{\vskip0.1em\noindent\hrulefill\vskip1em\begin{ExerciseList}}
\newcommand{\eexl}{\end{ExerciseList}\hrulefill}

\newcommand{\bthm}{\begin{theorem}}
\newcommand{\ethm}{\end{theorem}}
\newcommand{\bpro}{\begin{prop}}
\newcommand{\epro}{\end{prop}}
\newcommand{\bcor}{\begin{corollary}}
\newcommand{\ecor}{\end{corollary}}
\newcommand{\bcon}{\begin{conjecture}}
\newcommand{\econ}{\end{conjecture}}
\newcommand{\bp}{\begin{proof}}
\newcommand{\ep}{\end{proof}}
\newcommand{\blem}{\begin{lemma}}
\newcommand{\elem}{\end{lemma}}
\newcommand{\bn}{\begin{note}}
\newcommand{\en}{\end{note}}
\newcommand{\benum}{\begin{enumerate}}
\newcommand{\eenum}{\end{enumerate}}
\newcommand{\bed}{\begin{defn}}
\newcommand{\eed}{\end{defn}}
\newcommand{\brem}{\begin{remark}}
\newcommand{\erem}{\end{remark}}

\newcommand{\btik}{\begin{tikzpicture}\begin{axis}[scale=0.5,axis y line=center, axis x line=middle]}
\newcommand{\etik}{\end{axis}\end{tikzpicture}}

\let\into=\hookrightarrow
\let\mapsto=\longmapsto

\newcommand{\upperRomannumeral}[1]{\uppercase\expandafter{\romannumeral#1}}

 \usepackage{tikz-cd}
\usetikzlibrary{decorations.pathmorphing}
\usepackage{graphicx}
\usepackage{stmaryrd}
\usepackage{url}
\usepackage[authoryear,sort&compress,square]{natbib}
\let\cite=\citep
	\usepackage{fancyhdr}
	\usepackage[colorlinks,citecolor=blue,pagebackref,hypertexnames=false]{hyperref} \usepackage{zref-clever}
	\let\Cref=\zcref
\usepackage{titlesec}

\usepackage{tocloft}
\advance\cftsecnumwidth 0.5em\relax
\advance\cftsubsecindent 0.5em\relax
\advance\cftsubsecnumwidth 0.5em\relax

\newtoggle{draft}
\toggletrue{draft}
\newtheorem{theorem}[equation]{Theorem}      \newtheorem{theoremdef}[equation]{Theorem-Definition}
\newtheorem{lemma}[equation]{Lemma}          \newtheorem{corollary}[equation]{Corollary}  \newtheorem{proposition}[equation]{Proposition}
\newtheorem{prop}[equation]{Proposition}

\theoremstyle{definition}
\newtheorem{conj}[equation]{Conjecture}
\newtheorem{conjecture}[equation]{Conjecture}

\theoremstyle{definition}
\newtheorem{defn}[equation]{Definition}

\theoremstyle{remark}

\theoremsymbol{\ensuremath{\bullet}}

\theoremstyle{definition}
\newtheorem{remark}[equation]{Remark}
\theoremsymbol{\ensuremath{\bullet}}

\newcommand{\bthmdef}{\begin{theoremdef}}
\newcommand{\ethmdef}{\end{theoremdef}}

\zcsetup{cap} \AddToHook{env/theorem/begin}{\zcsetup{countertype={equation=theorem}}}
\AddToHook{env/lemma/begin}{\zcsetup{countertype={equation=lemma}}}
\AddToHook{env/proposition/begin}{\zcsetup{countertype={equation=proposition}}}
\AddToHook{env/prop/begin}{\zcsetup{countertype={equation=prop}}}
\AddToHook{env/theoremdef/begin}{\zcsetup{countertype={equation=theoremdef}}}
\AddToHook{env/defn/begin}{\zcsetup{countertype={equation=defn}}}
\AddToHook{env/definition/begin}{\zcsetup{countertype={equation=definition}}}
\AddToHook{env/remark/begin}{\zcsetup{countertype={equation=remark}}}
\AddToHook{env/rem/begin}{\zcsetup{countertype={equation=rem}}}
\AddToHook{env/conj/begin}{\zcsetup{countertype={equation=conj}}}
\AddToHook{env/conjecture/begin}{\zcsetup{countertype={equation=conjecture}}}
\AddToHook{env/question/begin}{\zcsetup{countertype={equation=question}}}
\AddToHook{env/example/begin}{\zcsetup{countertype={equation=example}}}

\zcRefTypeSetup{equation}{
	Name-sg = eq. ,
	name-sg = eq. ,
	Name-pl = eqns. ,
	name-pl = eqns. ,
}

\zcRefTypeSetup{theorem}{
	Name-sg = Theorem ,
	name-sg = theorem ,
	Name-pl = Theorems ,
	name-pl = theorems ,
}

\zcRefTypeSetup{theoremdef}{
	Name-sg = Theorem-Definition ,
	name-sg = theorem-definition ,
	Name-pl = Theorem-Definitions ,
	name-pl = theorem-definitions ,
}

\zcRefTypeSetup{proposition}{
	Name-sg = Proposition ,
	name-sg = proposition ,
	Name-pl = Propositions ,
	name-pl = propositions ,
}

\zcRefTypeSetup{lemma}{
	Name-sg = Lemma ,
	name-sg = lemma ,
	Name-pl = Lemmas ,
	name-pl = lemmas ,
}

\zcRefTypeSetup{remark}{
	Name-sg = Remark ,
	name-sg = remark ,
	Name-pl = Remarks ,
	name-pl = remarks ,
}

\zcRefTypeSetup{question}{
	Name-sg = Question ,
	name-sg = question ,
	Name-pl = Questions ,
	name-pl = questions ,
}

\zcRefTypeSetup{example}{
	Name-sg = Example ,
	name-sg = example ,
	Name-pl = Examples ,
	name-pl = examples ,
}

\titleformat{\subsection}[runin]{\normalfont\bfseries}{\S\ \thesubsection}{.5em}{}[{\ \ }]
\titlespacing{\subsection}{0pt}{1.5ex plus .1ex minus .2ex}{0pt}
\titleformat{\subsubsection}[runin]{\normalfont\bfseries}{\S\ \thesubsubsection}{.5em}{}[{\ \ }]
\titlespacing{\subsubsection}{0pt}{1.5ex plus .1ex minus .2ex}{0pt}

\let\into=\hookrightarrow
\let\isom=\simeq

\newcommand{\bQ}{{\bar{\Q}}}

\newcommand{\C}{{\mathbb C}}

\newcommand{\disc}{\mathop{{\rm disc}}}

\newcommand{\gal}{{\rm Gal}}
\newcommand{\Gm}{\mathbb{G}_m}

\newcommand{\Q}{{\mathbb Q}}
\newcommand{\R}{{\mathbb R}}

\newcommand{\Z}{{\mathbb Z}}

\renewcommand{\O}{{\mathcal O}}

\renewcommand{\wp}{{\mathfrak p}}

\newcommand{\invlim}{\varprojlim}
\newcommand{\mapright}[1]{{\xymatrix{{}\ar[r]^{#1}&{}}}}

\renewcommand{\bpro}{\begin{proposition}}
	\renewcommand{\epro}{\end{proposition}}
\newcommand{\bdefn}{\begin{defn}}
	\newcommand{\edefn}{\end{defn}}
\renewcommand{\bcon}{\begin{conj}}
	\renewcommand{\econ}{\end{conj}}

\title{The Categorical Local Langlands Correspondence and Anabelomorphy}
\author{Kirti Joshi}

\zcsetup{cap} \AddToHook{env/theorem/begin}{\zcsetup{countertype={equation=theorem}}}
\AddToHook{env/corollary/begin}{\zcsetup{countertype={equation=corollary}}}
\AddToHook{env/lemma/begin}{\zcsetup{countertype={equation=lemma}}}
\AddToHook{env/proposition/begin}{\zcsetup{countertype={equation=proposition}}}
\AddToHook{env/prop/begin}{\zcsetup{countertype={equation=prop}}}
\AddToHook{env/theoremdef/begin}{\zcsetup{countertype={equation=theoremdef}}}
\AddToHook{env/defn/begin}{\zcsetup{countertype={equation=defn}}}
\AddToHook{env/definition/begin}{\zcsetup{countertype={equation=definition}}}
\AddToHook{env/remark/begin}{\zcsetup{countertype={equation=remark}}}
\AddToHook{env/rem/begin}{\zcsetup{countertype={equation=rem}}}
\AddToHook{env/conj/begin}{\zcsetup{countertype={equation=conj}}}
\AddToHook{env/conjecture/begin}{\zcsetup{countertype={equation=conjecture}}}
\AddToHook{env/question/begin}{\zcsetup{countertype={equation=question}}}
\AddToHook{env/example/begin}{\zcsetup{countertype={equation=example}}}

\zcRefTypeSetup{equation}{
	Name-sg = eq. ,
	name-sg = eq. ,
	Name-pl = eqns. ,
	name-pl = eqns. ,
}

\zcRefTypeSetup{theorem}{
	Name-sg = Theorem ,
	name-sg = theorem ,
	Name-pl = Theorems ,
	name-pl = theorems ,
}

\zcRefTypeSetup{theoremdef}{
	Name-sg = Theorem-Definition ,
	name-sg = theorem-definition ,
	Name-pl = Theorem-Definitions ,
	name-pl = theorem-definitions ,
}

\zcRefTypeSetup{corollary}{
	Name-sg = Corollary ,
	name-sg = corollary ,
	Name-pl = Corollaries ,
	name-pl = corollaries ,
}

\zcRefTypeSetup{proposition}{
	Name-sg = Proposition ,
	name-sg = proposition ,
	Name-pl = Propositions ,
	name-pl = propositions ,
}

\zcRefTypeSetup{prop}{
	Name-sg = Proposition ,
	name-sg = proposition ,
	Name-pl = Propositions ,
	name-pl = propositions ,
}

\zcRefTypeSetup{lemma}{
	Name-sg = Lemma ,
	name-sg = lemma ,
	Name-pl = Lemmas ,
	name-pl = lemmas ,
}

\zcRefTypeSetup{remark}{
	Name-sg = Remark ,
	name-sg = remark ,
	Name-pl = Remarks ,
	name-pl = remarks ,
}

\zcRefTypeSetup{question}{
	Name-sg = Question ,
	name-sg = question ,
	Name-pl = Questions ,
	name-pl = questions ,
}

\zcRefTypeSetup{example}{
	Name-sg = Example ,
	name-sg = example ,
	Name-pl = Examples ,
	name-pl = examples ,
}

\zcRefTypeSetup{defn}{
	Name-sg = Definition ,
	name-sg = definition ,
	Name-pl = Definitions ,
	name-pl = definitions ,
}

\zcRefTypeSetup{conj}{
	Name-sg = Conjecture ,
	name-sg = conjecture ,
	Name-pl = Conjectures ,
	name-pl = conjectures ,
}

\zcRefTypeSetup{conjecture}{
	Name-sg = Conjecture ,
	name-sg = conjecture ,
	Name-pl = Conjectures ,
	name-pl = conjectures ,
}

\newcommand{\fixnumberwithin}[1]{
\numberwithin{equation}{#1}
	\numberwithin{theorem}{#1}
	\numberwithin{conj}{#1}
	\numberwithin{conjecture}{#1}
	\numberwithin{lemma}{#1}
	\numberwithin{proposition}{#1}
	\numberwithin{prop}{#1}
	\numberwithin{corollary}{#1}
	\numberwithin{defn}{#1}
	\numberwithin{definition}{#1}
	\numberwithin{remark}{#1}
	\numberwithin{rem}{#1}
	\numberwithin{question}{#1}
}

\newcommand{\nws}{\fixnumberwithin{section}}
\newcommand{\nwss}{\fixnumberwithin{subsection}}
 
\begin{document}
\maketitle
\nwss

\lhead{}

\newcommand{\benumlab}{\begin{enumerate}[label={{\bf(\arabic{*})}}]}

\newcommand{\bedef}{\begin{defn}}
\newcommand{\eedef}{\end{defn}}

 \newcommand{\bL}{\overline{L}}
 \newcommand{\bkm}{\bK_M}
 \newcommand{\vbk}{v_{\bK}}
 \newcommand{\vbkm}{v_{\bkm}}
\newcommand{\ocs}{\O^\circledast}
\newcommand{\ot}{\O^\triangleright}
\newcommand{\ocsk}{\ocs_K}
\newcommand{\otk}{\ot_K}
\newcommand{\ok}{\O_K}
\newcommand{\oko}{\O_K^1}
\newcommand{\oks}{\ok^*}
\newcommand{\Qpb}{\overline{\Q}_p}
\newcommand{\Qpbh}{\widehat{\overline{\Q}}_p}
\newcommand{\tr}{\triangleright}
\newcommand{\ocpt}{\O_{\C_p}^\tr}
\newcommand{\ocpf}{\O_{\C_p}^\flat}
\newcommand{\sG}{\mathscr{G}}
\newcommand{\sY}{\mathscr{Y}}
\newcommand{\sxqp}{\mathscr{X}_{\cpt,\Q_p}}
\newcommand{\syqp}{\mathscr{Y}_{\cpt,\Q_p}}
\newcommand{\sxfe}{\mathscr{X}_{F,E}}
\newcommand{\sxfep}{\mathscr{X}_{F,E'}}
\newcommand{\syfe}{\mathscr{Y}_{F,E}}
\newcommand{\syfep}{\mathscr{Y}_{F,E'}}
\newcommand{\loglt}{\log_{\sG}}
\newcommand{\fc}{\mathfrak{t}}
\newcommand{\ku}{K_u}
\newcommand{\kup}{\ku'}
\newcommand{\kt}{\tilde{K}}
\newcommand{\sGpf}{\sG(\O_K)^{pf}}
\newcommand{\hgm}{\widehat{\mathbb{G}}_m}
\newcommand{\hgmp}{\widehat{\mathbb{G}}_{m;p}}
\newcommand{\bE}{\overline{E}}

\newcommand{\sM}{\mathscr{M}}

\newcommand{\vlnon}{\mathbb{V}_L^{non}}
\newcommand{\vlarch}{\mathbb{V}_L^{arc}}

\newcommand{\fA}{\mathfrak{A}}

\newcommand{\V}{\mathbb{V}}
\newcommand{\vl}{\V_L}
\newcommand{\lvbh}{\widehat{\overline{L}}_v}
\newcommand{\lvbht}{{\widehat{\overline{L}}_v}^\flat}
\newcommand{\lvbhmax}{\widehat{\overline{L}}_v^{max}}
\newcommand{\lvbhtmax}{{\widehat{\overline{L}}_v}^{max\flat}}
\newcommand{\lvbhreal}{\widehat{\overline{L}}_v^{\R}}
\newcommand{\lvbhtreal}{{\widehat{\overline{L}}_v}^{\R\flat}}

\newcommand{\sX}{\mathscr{X}}
\newcommand{\syflv}{\sY_{\lvbht,L_v}}
\newcommand{\syflvmax}{\sY_{\lvbhtmax,L_v}}
\newcommand{\sxflv}{\sX_{\lvbht,L_v}}
\newcommand{\sxflvmax}{\sX_{\lvbhtmax,L_v}}
\newcommand{\bsY}{\mathscrbf{Y}}
\newcommand{\bsX}{\mathscrbf{X}}
\newcommand{\yadl}{\bsY_L}
\newcommand{\yadlmax}{\bsY_L^{max}}
\newcommand{\yadq}{\bsY_{\Q}}
\newcommand{\yadqmax}{\bsY_{\Q}^{max}}
\newcommand{\xadl}{\bsX_L}
\newcommand{\xadq}{\bsX_{\Q}}
\newcommand{\xadlmax}{\bsX_L^{max}}
\newcommand{\xadqmax}{\bsX_{\Q}^{max}}
\newcommand{\yadlpoint}{\{(L_\wp\into K_{\wp}, K_\wp\isom \cpt)\}_{\wp\in\vlnon}}
\newcommand{\xadlpoint}{\{(L_\wp\into K_{\wp}, K_\wp\isom \cpt)\}_{\wp\in\vlnon}}

\newcommand{\bv}{\bar{v}}
\newcommand{\sGal}[1]{\mathbf{G}_{#1}}

\newcommand{\arith}[1]{\mathfrak{arith}(#1)}
\newcommand{\arithl}{\arith{L}}
\newcommand{\adel}[1]{\mathfrak{adel}(#1)}
\newcommand{\adell}{\adel{L}}
\newcommand{\by}{{\bf y}}

\newcommand{\fet}{\mathcal{F\hat{e}t}}
\newcommand{\sI}{\mathcal{I}}
\newcommand{\sA}{\mathcal{A}}

\newcommand{\logp}{\log^+}
\newcommand{\sFrob}{\mathcal{F\!r\!o\!b}}
\newcommand{\sF}{\mathscr{F}}
\newcommand{\pow}[2]{#1\llbracket#2\rrbracket}

\newcommand{\ainf}{A_{\textrm{inf}}}

\newcommand{\ells}{{\ell^*}}
\newcommand{\glqp}{{\rm GL}_2^+(\Q)}
\newcommand{\tslq}{\widetilde{{\rm SL}_2(\Q)}}
\newcommand{\tslql}{{\widetilde{{\rm SL}_2(\Q)}}^\ells}
\newcommand{\slq}{{\rm SL_2}(\Q)}
\newcommand{\tglqp}{\widetilde{{\rm GL}_2^+(\Q)}}
\newcommand{\tglqpl}{{\widetilde{{\rm GL}_2^+(\Q)}}^\ells}
\newcommand{\slr}{{\rm SL}_2(\R)}
\newcommand{\tslr}{\widetilde{{\rm SL}_2(\R)}}
\newcommand{\slz}{{\rm SL}_2(\Z)}
\newcommand{\tslz}{\widetilde{{\rm SL}_2(\Z)}}
\newcommand{\fH}{\mathfrak{H}}
\newcommand{\syi}{\sY_{\infty}}
\newcommand{\syis}{\sY_{\infty,s}}
\newcommand{\syils}{\sY_{\infty,s}^\ells}
\newcommand{\sxi}{\sX_\infty}
\newcommand{\vphi}{\varphi}
\newcommand{\flog}{\mathfrak{log}}
\newcommand{\vnl}{\V_L^{non}}
\newcommand{\Bl}{\mathbb{B}_L}
\newcommand{\Blp}{\mathbb{B}^+_L}
\newcommand{\val}{\V_L^{arc}}
\newcommand{\ssep}{\S\,}
\newcommand{\bvphi}{\boldsymbol{\varphi}}

\newcommand{\sL}{\mathscr{L}}
\newcommand{\sE}{\mathscr{E}}
\newcommand{\sT}{\mathscr{T}}
\newcommand{\onto}{\twoheadrightarrow}

\togglefalse{draft}

\iftoggle{draft}{\pagewiselinenumbers}{\relax}

\begin{abstract}
Let $G/\Q_p$ be a  connected, split, reductive group over $\Q_p$.
In this paper I show that if $K$ and $L$ are anabelomorphic $p$-adic fields i.e. $K$ and $L$ have topologically isomorphic absolute Galois groups, then the stacks of Langlands parameters (for the fields $K$ and $L$) considered in \cite{fargues-scholze}, are also isomorphic (\Cref{thm:anab1}). This leads to \Cref{con:synchronization-conj} which provides a precise relationship between the main conjecture of \cite{fargues-scholze} and anabelomorphy of $p$-adic fields considered in \cite{joshi-anabelomorphy}. I establish my conjecture for a split torus in \Cref{thm:anab2}.
\end{abstract}

\setcounter{tocdepth}{3}
\tableofcontents

\newcommand{\bK}{\bar{K}}
\renewcommand{\bL}{\bar{L}}
\newcommand{\gamk}{\Gamma_K}
\newcommand{\gaml}{\Gamma_L}
\newcommand{\bZ}{\bar{\Z}}
\newcommand{\anab}{\leftrightsquigarrow}
\newcommand{\anabmapright}[1]{\stackrel{#1}{\leftrightsquigarrow}}
\newcommand{\gd}{\hat{G}}
\newcommand{\td}{\hat{T}}
\newcommand{\pargk}{\text{Par}_{G;K}}
\newcommand{\pargl}{\text{Par}_{G;L}}
\newcommand{\indcoh}{\text{IndCoh}}

\newcommand{\sD}{\mathcal{D}}
\newcommand{\bungk}{\text{Bun}_{G;K}}
\newcommand{\bungl}{\text{Bun}_{G;L}}
\newcommand{\buntk}{\text{Bun}_{T;K}}
\newcommand{\buntl}{\text{Bun}_{T;L}}

\newcommand{\dbungk}{\sD_{lis}(\text{Bun}_{G;K},\bQ_\ell)}
\newcommand{\dbungl}{\sD_{lis}(\text{Bun}_{G;L},\bQ_\ell)}

\newcommand{\dbuntk}{\sD_{lis}(\text{Bun}_{T;K},\bQ_\ell)}
\newcommand{\dbuntl}{\sD_{lis}(\text{Bun}_{T;L},\bQ_\ell)}
\newpage

\nwss
\section{Introduction}\label{se:intro}
The Local Langlands Correspondence in all its avatars (\cite{henniart1988}, \cite{fargues-scholze}, \citep*{emerton2025}) deals with representations of the absolute Galois groups (or more precisely of the Weil or the Weil-Deligne groups, of $p$-adic fields) and (suitably defined) automorphic representations. On the other hand, there exists non-isomorphic $p$-adic fields $K,L$ with topologically isomorphic absolute Galois groups. This leads to the notion of \emph{anabelomorphy}, studied in detail in \cite{joshi-anabelomorphy}. More precisely, two $p$-adic fields $K,L$ with topologically isomorphic absolute Galois groups are said to be \emph{anabelomorphic $p$-adic fields}  (examples of such fields are given by \cite[Theorem 3.5.1 and Lemma 4.4]{joshi-anabelomorphy}).  Especially, from  \cite[Proposition 7.1,2]{joshi-anabelomorphy}, one sees that anabelomorphic $p$-adic fields have isomorphic Weil (and Weil-Deligne) groups. In particular, anabelomorphic $p$-adic fields have naturally equivalent categories of finite dimensional Galois representations. 

So one may ask if the automorphic sides of the Local Langlands correspondence for anabelomorphic fields are also directly and naturally related.  This question was posed in \cite[\ssep 7]{joshi-anabelomorphy}, where  partial, affirmative answers were obtained in the context of \cite{henniart1988} for all principal series representations and all supercuspidal representations in the tame case ($p\nmid n$) for $GL_n$ (see \cite[Theorem 7.3.3 and Theorem 7.4.1]{joshi-anabelomorphy}). In \cite[Theorem 7.6.1]{joshi-anabelomorphy}, one obtains a complete result for $GL_2$ and all odd residue characteristics. For my discussion of the $p$-adic context of \cite{emerton2025}, see \cite[\ssep 7.7, and also Question 15.17]{joshi-anabelomorphy}.

In the present paper, I formulate the precise expectation (\Cref{con:synchronization-conj}) for the categorical avatar given by \cite{fargues-scholze} for connected, split reductive groups over $\Q_p$ considered over anabelomorphic $p$-adic fields. In \Cref{thm:anab1} (also see \Cref{cor:anab-2}), I establish the equivalence of the Galois side i.e. equivalence of the stacks of Langlands parameters considered in \cite{fargues-scholze}. In \Cref{pr:schematic-isom}, I establish the isomorphism of schemes of Langlands parameters defined by \cite[Theorem VIII.1.3]{fargues-scholze}.  Equivalence of Whittaker data is established in \Cref{pr:whittaker}. Finally in \Cref{thm:anab2}, \Cref{con:synchronization-conj} is proved for split torii.

\emph{The Asides \ssep3.4--\ssep3.7 may be skipped on initial reading: \ssep3.4 is about groups with no non-trivial pure inner forms; \ssep3.5 clarifies how anabelomorphy plays out on the $G$-bundle side of \cite{fargues-scholze} for extended inner forms of $G=GL_n$; \ssep3.6 (resp. \ssep3.7) deals with anabelomorphy and quadratic (resp. Hermitian) spaces.}

\subsection{Acknowledgements} I thank Peter Scholze and Konrad Zou for providing comments, suggestions and references.

\section{The main theorem in the split case}
\nwss
\subsection{Anabelomorphy of $p$-adic fields}
Let $K,L$ be $p$-adic fields (i.e. $K,L$ are finite extensions of $\Q_p$ for some prime number $p$), and  choose algebraic closures $\bK\supset K$ and $\bL\supset L$ of $K$ and $L$. Let $\gamk=\gal(\bK/K)$ (resp. $\gaml=\gal(\bL/L)$) be their absolute Galois groups computed using  $\bK$ (resp. $\bL$).  Following \cite[Definition 2.1.1]{joshi-anabelomorphy}, one says that $K,L$ are \emph{anabelomorphic} if there exists an isomorphism of topological groups $\alpha:\gamk\mapright{\isom} \gaml$, and refer to $\alpha$ as an anabelomorphism between $K$ and $L$.  An anabelomorphism will be denoted as  $K\anabmapright{\alpha} L$ or if $\alpha$ is not important, simply as $K\anab L$. 
Notably, anabelomorphy of $p$-adic field is an equivalence relation. For basic properties of anabelomorphisms, see \cite[\ssep 2]{joshi-anabelomorphy}; for examples of anabelomorphic, but non-isomorphic, $p$-adic fields see \cite[Lemma 4.4]{joshi-anabelomorphy} and \cite[Theorem 3.5.1]{joshi-anabelomorphy}. 

Following \cite[Definition 2.1.3]{joshi-anabelomorphy}, a quantity, property, or an algebraic structure associated with $K$ is said to be \emph{amphoric} if it is common to all members of the anabelomorphism class of $K$. For examples of amphoric quantities, properties and structures, see \cite[Theorem 3.4.1]{joshi-anabelomorphy} and especially \cite[Theorem 3.6.1]{joshi-anabelomorphy}. Notably, the prime number $p$ and cardinality $q$ of the residue field $K$ is amphoric (\cite[Theorem 3.4.1 and Proposition 7.1.2]{joshi-anabelomorphy}). 

Suppose that $K\anab L$ are anabelomorphic $p$-adic fields. By \cite[Theorem 3.5.1]{joshi-anabelomorphy}, one sees that the maximal abelian subfields $K\supseteq K^0\supset \Q_p$ (resp. $L\supseteq L^0\supset \Q_p$) contained in $K$ (resp. $L$) are isomorphic. 

In many proofs involving anabelomorphy, one encounters or constructs natural bijections, or natural equivalences, between various sets or categories associated with anabelomorphic fields $K$, $L$ (in particular, these sets and categories are amphoric). I will habitually use the terms \emph{synchronize} or \emph{synchronization} for such constructions. [The term `synchronization' is short for `anabelomorphic synchronization.'] An example for the present paper is the synchronization of the Whittaker data (\Cref{pr:whittaker}) and the Synchronization Conjecture~\ref{con:synchronization-conj}. Other examples of synchronizations can be found in \cite{joshi-anabelomorphy}.

\subsection{Anabelomorphy and the stacks of Langlands Parameters} 
Now fix a prime number $\ell\neq p$, and let $\Z_\ell,\Q_\ell$ have their standard meaning. Let $\bQ_\ell\supset \Q_\ell$ be an algebraic closure, and $\bZ_\ell\supset\Z_\ell$ be the integral closure of $\Z_\ell$ in $\bQ_\ell$. Fix a $\sqrt{q}\in\bQ_\ell$. The notion of an \emph{animated ring} is as defined in \cite[5.1.4]{scholze2023-cesnavicius}. The next theorem establishes a natural isomorphism between stacks of Langlands parameters:
\bthm\label{thm:anab1} 
Let $L\anabmapright{\alpha} K$ be an  anabelomorphism  of $p$-adic fields.  Let $G$ be a connected, split, reductive algebraic group over $\Q_p$ (one can assume that $G$ has no non-trivial pure inner forms e.g. \Cref{pr:no-inner-form}). Let $\gd$ be the  dual group of $G$. Then $\alpha$ induces an isomorphism of the $\bZ_\ell$-stacks of Langlands parameters
\be \alpha:\pargl \mapright{\isom} \pargk.\ee
In other words, $\pargk$ is amphoric.
\ethm
\bcor\label{cor:anab-2}
In the notation and assumptions of \Cref{thm:anab1}. The anabelomorphism $L\anabmapright{\alpha} K$ induces a natural equivalence of categories
\be  
\indcoh(\pargl)\mapright{\isom} \indcoh(\pargk).
\ee
In other words, $\indcoh(\pargl)$ is also amphoric.
\ecor
\bp[Proof of \Cref{cor:anab-2}]
This is immediate from \Cref{thm:anab1} as the two stacks are isomorphic and the general formalism of ind-coherent sheaves given in \cite{gaitsgory2017}.
\ep

\bp[Proof of \Cref{thm:anab1}]
For the construction of stack of Langlands parameters see \citep*{zhu2020}, \citep*{dat2020},  \cite[Chapter VIII]{fargues-scholze}, \cite[3.2]{hansen2026} (I will follow the latter two).

By \cite[Theorem 1.3]{conrad2015-non-split}, any connected, split, reductive group over $\Q_p$ arises, by base-change, from a Chevalley group scheme over $\Z$, which is unique up to a $\Z$-isomorphism. Note that the formation of the dual $\gd$ is fixed by considering any algebraically closed field with a field embedding of $K$ and is independent of the topology of $K$. Hence, one can refer to $\gd$ independently of the choice of algebraic closures $\bK$ and $\bL$ of $K$ and $L$ respectively.  The existence of the Chevalley group scheme allows one to mitigate the difficulty highlighted in \Cref{re:abelianized-cohom} in the context of $\gd$ in the considerations which follow.

Since $G$ is split over $\Q_p$,  so is $G_L$ (resp. $G_K$). Strictly speaking, one should write $G_L=G\times_{\Q_p}L$ (resp. $G_K=G\times_{\Q_p}K$). But for notational simplicity, the subscripts will be suppressed. Since $G$ is split over $L$ (resp. $K$) means that the action of $W_L$ on $\gd_L$ (resp. action of $W_K$ on $\gd_K$) is trivial.

Let $\gamk \supset I_K \supset P_K$ (resp. $\gaml \supset I_L \supset P_L$) be the inertia and the wild inertia subgroups of $\gamk$ (resp. $\gaml$). Let $W_K\subset \gamk$ (resp. $W_L\subset \gaml$) be the Weil group of $K$ (resp. $L$). 

The key point in the proof is that by \cite[Proposition 7.1.2]{joshi-anabelomorphy} Weil group is amphoric i.e. for any  anabelomorphism $L\anabmapright{\alpha}K$ induces a topological isomorphism of Weil groups 
\be 
\begin{tikzcd}
	W_L \ar[r,bend left,"\alpha","\isom"'] & W_K \ar[l,bend left,"\alpha^{-1}","\isom"'].
\end{tikzcd}
\ee
and by \cite[Theorm 3.6.1]{joshi-anabelomorphy}, $\alpha$ maps $I_K$ isomorphically onto $I_L$ and $P_K$ isomorphically onto $P_L$ and similarly for $\alpha^{-1}$. 

From definition of $\pargk$ given in \cite[Chapter VIII.1.1]{fargues-scholze},  one knows that $\pargk$ depends only on the Weil group (considered as a condensed animated group).  So proof of the theorem will be self-evident for the experts from the amphoricty of the Weil group. But let me walk the readers through the details for added clarity.

To prove the equivalence between 
$\pargl$ and $\pargk$, it will be convenient to work with Langlands parameters for both the sides. One starts with a Langlands parameter for $K$  and from it one constructs a Langlands parameter for $L$  directly and vice versa. 

Let $A$ be any animated ring over $\bZ_\ell$ in the sense of \cite{scholze2023-cesnavicius}. Then a Langlands parameter for $K$ is a homomorphism of condensed animated groups
\be 
\phi:W_K\to \gd(A).
\ee
[Note that this is where one uses the assumption that $G$ is split, otherwise one works with 1-cocycles \cite[Definition VIII.1.1]{fargues-scholze}.] 

By \cite{joshi-anabelomorphy}, $\alpha$ induces an isomorphism $W_L\mapright{\isom} W_K$. Hence precomposing $\phi$ with $\alpha$ gives us a homomorphism
\be 
\phi':W_L\mapright{\isom} W_K\mapright{\phi} \gd(A).
\ee
Conversely, if $\phi':W_L\to \gd(A)$ is a Langlands parameter for $L$, then one constructs a Langlands parameter $\phi$ by precomposition with the isomorphim $\alpha^{-1}:W_K\mapright{\isom}W_L$. 
Thus one has a natural isomorphism between 
\be\label{eq:phi-phiprime}
Hom(W_K,\gd(A))\mapright{\phi\mapsto \phi'} Hom(W_L,\gd(A)).
\ee

The stack $\pargk$ (resp. $\pargl$) is the stack in infinity groupoids
\be 
A\mapsto Hom(W_K,\gd(A))/\gd(A)
\ee
where the action of $\gd(A)$ is by conjugation of elements in $\gd(A)$. Thus \eqref{eq:phi-phiprime} leads to the natural bijection
\be  
Hom(W_K,\gd(A))/\gd(A)\mapright{\alpha^{-1}}Hom(W_L,\gd(A))/\gd(A)
\ee
given by \be\phi\bmod{\gd(A)}\mapsto \phi'\bmod{\gd(A)}.\ee
Hence, the asserted equivalence $\alpha:\pargl\mapright{\isom}\pargk$ between the stacks  of Langlands parameters is established. This proves the theorem. 
\ep

The following result leads to a more refined version of \Cref{thm:anab1}.

\bpro\label{pr:schematic-isom}
Let $L\anabmapright{\alpha}K$ be an anabelomorphism of $p$-adic fields and let $G$  satisfy all the hypothesis of \Cref{thm:anab1}. Then $\alpha$ induces a natural isomorphism of affine $\Z_\ell$-schemes of Langlands parameters
\be 
Z^1(W_K,\gd)\mapright{\alpha^{-1}} Z^1(W_L,\gd).
\ee
\epro
\bp 
Experts will recognize that this assertion is also a consequence of the amphoricity of the Weil group, but let me give a detailed proof for completeness. The explicit construction of $Z^1(W_K,\gd)$ and its properties is given by \cite[Theorem VIII.1.3]{fargues-scholze}. That construction also shows that $Z^1(W_K,\gd)$ is a union of open and closed affine $\Z_\ell$-schemes $Z^1(W_K/P',\gd)$ where $P'\subset P_K\subset I_K\subset W_K$ runs through the collection of open subgroups of the wild inertia subgroup $P_K$.

By \cite[Theorem 3.6.1]{joshi-anabelomorphy}, the inertia subgroup and the wild inertia subgroups $P_K\subset I_K$ of $W_K$ are amphoric. Thus, as $P'\subset P_K$ runs through open subgroups of $P_K$, $\alpha^{-1}(P)\subseteq \alpha^{-1}(P_K)=P_L$ runs through the open subgroups of $P_L$. Thus one has a natural isomorphism of $\Z_\ell$-schemes 
\be Z^1(W_K/P',\gd) \isom Z^1(W_L/\alpha^{-1}(P'),\gd),\ee
which induces an isomorphism of their disjoint unions
\be Z^1(W_K,\gd) \isom Z^1(W_L,\gd).\ee
This completes the proof.
\ep

\section{The synchronization conjecture}
\nwss
\subsection{The conjecture of \cite{fargues-scholze}}
For any $p$-adic field $K$ of residue cardinality $q$ such that $\sqrt{q}\in\bQ_\ell$, any quasi-split reductive group $G$ over $K$, and any generic character $\psi:U(K)\to\bQ_\ell^*$ of the unipotent radical $U\subset B\subset G$ in \cite{fargues-scholze}, it is conjectured that there exists an equivalence, with certain natural properties, between the categories 
\be
\begin{tikzcd}
		\dbungk \ar[r,"FS_{\psi;K}","\isom"'] & \indcoh(\pargk).
\end{tikzcd} 
\ee
The equivalence $FS_{\psi;K}$ is denoted by $L_\psi$ in \cite{fargues-scholze} (and their notation is avoided here as $L$ is a $p$-adic field throughout this paper, and instead the   notation used here mostly follows \cite{hansen2026}). I will not define $\dbungk$ here as it is quite elaborate and readers are referred to \cite[1.5]{fargues-scholze} or \cite[4.2]{hansen2026}; the definition of Ind-coherent sheaves on stacks is also quite elaborate and the reader is referred to \cite{gaitsgory2017}.

This is the main conjectural version of the categorical local Langlands Correspondence of \cite{fargues-scholze}. In \cite{zouKonrad2024}, this settled for torii and in \cite{hansen2026}, this conjectured settled under some reasonable hypothesis for many split groups including $GL_n$. 

\subsection{Synchronization of the Whittaker data}
The next proposition deals with anabelomorphic synchronization of Whittaker datum of \cite[Conjecture I.10.2 and Chapter X]{fargues-scholze}. Recall from \cite{borel-linear-algebraic} or \cite{humphreys-linear-algebraic} the following. Let $G$ be a connected, split reductive group over $\Q_p$, $B\subset G$ be a Borel subgroup. Let $U\subseteq B\subseteq G$ be the unipotent radical of $B$.  Write $G_K=G\times_{\Q_p}K$ and similarly for $G_L$. Then as $K,L$ are perfect fields, $B_K$ (resp. $B_L$) is a Borel subgroup of $G_K$ (resp. $G_L$) and $U_K$ (resp. $U_L$) is the unipotent radical of $B_K$ (resp. $B_L$).

\newcommand{\Ga}{\mathbb{G}_a}
\bpro\label{pr:whittaker} 
Let $L\anabmapright{\alpha}K$ be an anabelomorphism of $p$-adic fields. Keep the above assumptions for $G$ and notations for the subgroups $U\subset B\subset G$ and their $K$ (resp. $L$) versions. Then $\alpha$ induces a natural bijection between the set of characters 
\be \psi:U_K(K)\to \bQ_\ell^*\ee and the set of characters 
\be\psi':U_L(L)\to\bQ_\ell^*.\ee
\epro
\bp 
Any homomorphism $U_K(K)\to \bQ_\ell^*$ factors through the abelianization of $U_K(K)$.
As $U$ is unipotent, $U/[U,U]=\Ga^m$ for some integer $m\geq 1$, and the exact sequence
\be 0\to [U_K,U_K]\to U_K\to \Ga^m\to 0,\ee
Taking Galois cohomology and noting that for any unipotent group $H^1(K,U_K)=0$ by \cite[Chapter 2, Lemma 2.7]{platonov-book}, one obtains
the exact sequence
\be 0\to [U_K(K),U_K(K)]\to U_K(K)\to \Ga^m(K)\to 0.\ee
and a similar exact sequence for $U_L(L)$. By \cite[Theorem 3.4.1]{joshi-anabelomorphy}, one has an isomorphism of  (additive) topological groups
\be 
U_L(L)/[U_L(L),U_L(L)]=L^m=\Ga^m(L)\isom \Ga^m(K)=K^m=U_K(K)/[U_K(K),U_K(K)].
\ee
Hence, one has a natural bijection of the sets of  characters $\psi: U_K(K)/[U_K(K),U_K(K)]\to\bQ_\ell^*$ and $\psi': U_L(L)/[U_L(L),U_L(L)]\to \bQ_\ell^*$
given using the above homeomorphism. This proves the assertion.
\ep
\newcommand{\sW}{\mathcal{W}}
\newcommand{\cInd}[2]{\text{c-Ind}_{#1}^{#2}}
\brem 
The datum $(U_K(K),\psi)$ is the Whittaker datum for $G_K$ (see \cite[Conjecture I.10.2 and Chapter X]{fargues-scholze}), and it gives rise to the Whitataker sheaf $\sW_K(U,\psi)$ on  $\bungk$ i.e. $\sW_K(U,\psi)$ is the sheaf corresponding to the compact induction $\cInd{U_K(K)}{G(K)}(\psi)$. 
\erem

\subsection{The Synchronization Conjecture}
In \cite{joshi-anabelomorphy}, I established a number of  synchronization results which relate the sets of irreducible admissible representation of $GL_n(K)$ and $GL_n(L)$. Specifically, \cite[Theorem 7.3.3]{joshi-anabelomorphy} shows that one can synchronize the sets of principal series representations of $GL_n(K)$ and $GL_n(L)$ for any pair of anabelomorphic $p$-adic fields. \cite[Theorem 7.4.1]{joshi-anabelomorphy} establishes the synchronization of supercuspidal representations (for these groups) assuming that $p\nmid n$ (i.e. in the tame case). Notably this allows one  to establish complete results for $GL_2$ for $K$ (and hence $L$) of odd residue characteristics (\cite[Theorem 7.6.1]{joshi-anabelomorphy}). That work was the motivation for \Cref{thm:anab1}, and it leads to the following:

\bcon[The synchronization conjecture]\label{con:synchronization-conj} Let $\ell\neq p$ be primes. Fix an algebraic closure $\bQ_\ell\supset \Q_\ell$ and let $\bZ_\ell\supset \Z_\ell$ be the integral closure of $\Z_\ell$ in $\bQ_\ell$.
Let $G/\Q_p$ be a connected, split reductive group. Let $K,L$ be anabelomorphic $p$-adic fields, and fix a $\sqrt{q}\in\bQ_\ell$. Then any anabelomorphism  $L\anabmapright{\alpha}K$ induces an equivalence of categories
\be 
\begin{tikzcd}
	\dbungl \ar[r,"\alpha_{Bun}","\isom"'] & \dbungk
\end{tikzcd}
\ee
such that one has a  natural commutative diagram of equivalences of categories
\be 
\begin{tikzcd}
	\dbungl \ar[d,"\alpha_{Bun}","\isom"']\ar[r,"FS_{\psi;L}","\isom"'] & \indcoh(\pargl)\ar[d, "\alpha_{Par}", "\isom"'] \\
		\dbungk \ar[r,"FS_{\psi;K}","\isom"'] & \indcoh(\pargk)
\end{tikzcd}
\ee
where $FS_{\psi;L}$ (resp. $FS_{\psi;K}$) is the conjectural categorical Langlands equivalence given by \cite[Conjecture X.1.4]{fargues-scholze} for $L$ (resp. $K$), and where $\alpha_{Par}$ is the equivalence given by \Cref{thm:anab1}. In other words, $\dbungk$ is amphoric. Moreover, $\alpha_{Bun}$ carries the Whittaker sheaf datum for $L$ to the Whittaker sheaf datum for $K$: 
\be 
\alpha_{Bun}(\sW_L(U_L,\psi'))=\sW_K(U_K,\psi)
\ee
via the synchronization of the Whittaker data given by \Cref{pr:whittaker}.
\econ

\brem\ 
\benumlab
\item In \cite[Proposition 1.4.5]{hansen2026} it is shown that there is at most one equivalence $FS_{\psi,K}$ whose existence is conjectured in \cite{fargues-scholze}. 
\item I do not know how to construct $\alpha_{Bun}$ for any pair of anabelomorphic $p$-adic fields beyond what is proved below.
\item The above statement will require dealing with pure inner forms. That needs a separate treatment from the point of view of anabelomorphy which is not discussed here. But for a large class of groups one can circumvent this issue via \Cref{pr:no-inner-form}.
\item In \cite{fargues-scholze} $\bungk$ is constructed via the construction of Fargues-Fontaine curves (\cite{fargues-fontaine}). On the other hand,  I have also shown that anabelomorphic $p$-adic fields lead to anabelomorphic Fargues-Fontaine curves (i.e. these curves have isomorphic \'etale fundamental groups), but these curves need not be isomorphic as $\Z$-schemes (see \cite{joshi-gconj}). So the existence of $\alpha_{Bun}$ seems quite subtle.
\eenum
\erem

\subsection{Aside: Groups with no non-trivial pure inner forms}
The following technical condition simplifies some of the proofs. In general, I do not expect that this is required.

\bpro\label{pr:no-inner-form}
Let $G$ be a split reductive group over a $p$-adic field $K$. Let $\bK\supset K$ be an algebraic closure of $K$ and let $\gamk=\gal(\bK/K)$  be the absolute Galois group of $K$.  
\benumlab
\item 
If $H^1(\gamk, G(\bK))={1}$, then $G$ has no non-trivial pure inner forms.
\item If $G$ is a finite direct product of any of the groups $$\left\{GL_n,SL_n, GSp_{2n}, Sp_{2n}, GSpin_{n}, Spin_{n}, E_8, F_4, G_2\right\},$$ then $G$ has no non-trivial pure inner forms.
\eenum
\epro
\bp 
By the definition of pure inner forms of $G$ given by\cite[Page 250]{conrad2014-reductive-book} it follows that if $H^1(\gamk,G(\bK))={1}$, then $G$ has no non-trivial pure inner forms. Hence, the assertion {\bf(1)} follows from the hypothesis. Now {\bf(1)}$\implies${\bf(2)} for all the groups in the list. To see this, one notes that the derived subgroup is the corresponding group without the letter $G$ in front of it (eg. $GL_n'=SL_n$ etc.). By \cite[14.2]{borel-linear-algebraic} the derived groups of all the groups in the list are semisimple. The remaining groups $E_8,F_4,G_2$ are simple and are their own derived groups. By \cite[8.1.11]{springer-linear-algebraic}, the simply connectedness of the derived groups of all the groups in the list can be verified by comparing root and weight lattices and this can be done case by case. 

Now the Galois cohomology sequence applied to 
\be 1\to G'\to G\to G/G'\to 1\ee
reduces the asserted vanishing for $G$ to that of its derived group $H^1(\gamk,G')=1$ and this assertion is \cite[Chapter 6, Theorem 6.4]{platonov-book}. 
\ep

\brem\label{re:abelianized-cohom} 
The above technical condition will necessary (at the moment) in some of the considerations below because of the following. The field $\bK$ together with its continuous $\gamk$-action is not amphoric (this is a consequence of \cite[Proposition 2.2]{mochizuki-local-gro} which implies that the completion $\widehat{\bK}$ equipped with its continuous $\gamk$-action is not amphoric in general). So it is generally not possible to choose an isomorphism of fields $\bK\isom \bL$ compatible with any anabelomorphism $K\anabmapright{\alpha}L$.  This leads to the following difficulty:  one cannot generally choose an isomorphism $G(\bK)\isom G(\bL)$ compatibly with anabelomorphisms and  the respective $\gamk,\gaml$ actions. In particular, one cannot relate the Galois cohomology groups $H^j(\gamk,G_K)$, $H^j(\gaml,G_L)$. On the other hand, in \cite[Corollary 5.4.1]{borovoi1998}, $H^1(\gamk,G_K)$ is identified with Borovoi's abelianized Galois cohomology $H^1_{ab}(\gamk,G_K)$ (and this group depends on the root datum via \cite[Corollary 5.5]{borovoi1998}) and hence one expects that $H^1_{ab}(\gamk,G_K)$ (and hence $H^1(\gamk,G_K)$) is amphoric.
\erem

\subsection{Aside: Synchronization of extended inner forms of $GL_n$}
The reductive group $GL_n$ has no pure inner forms by \Cref{pr:no-inner-form}. On the other hand $GL_n$ has extended inner forms $GL_m(D)$ for suitable $m$ and suitable central division algebras $D$. The following provides synchronization of extended inner forms of $GL_n$:

\bpro\label{pr:extended-syn} 
Any anabelomorphism $K\anabmapright{\alpha}L$ provides a synchronization of extended inner forms of $GL_n/K$ and $GL_n/L$ which takes a central division algebra $D/K$ with invariant $\mu\in\text{Br}(K)$ to a central division algebra $D'/L$ with the same invariant in $\text{Br}(L)$ via the amphoricity of Brauer groups given by \cite[Proposition 7.5.2]{joshi-anabelomorphy}.  This provides the synchronization of extended inner forms of $GL_n/K$, it maps
$GL_m(D)$ of $GL_n/K$ to the extended inner form $GL_m(D')$ of $GL_n/L$.
\epro 
\bp 
This is immediate from the amphoricity of the Brauer group given by \cite[Proposition 7.5.2]{joshi-anabelomorphy}. For completeness, let me prove the amphoricity of Brauer groups:
\be 
\text{Br}(K)=H^2(\gamk,\bK^*)\isom H^2(\gaml,\bL^*)=\text{Br}(L)
\ee
which follows from the amphoricity of the triple $\gamk\curvearrowright\bK^*$ given by \cite[Proposition 4.2]{hoshi-mono}.
\ep

\subsection{Aside: Synchronization of quadratic spaces}
\bpro 
Any anabelomorphism $K\anabmapright{\alpha}L$ provides a synchronization of non-singular or regular quadratic spaces over $K$ and $L$ respectively.
\epro
\bp 
By \cite[Chapter 5, 6.4 Corollary]{scharlau2012-quadratic-book}, a non-singular or regular quadratic space $(V,Q)$ over $K$ is determined, up to isomorphism, by three invariants $\dim(V)$, $\text{disc}(Q)$ and the Hasse invariant $\varepsilon(V,Q)$. One can assume that $Q=\langle a_1,\ldots,a_n\rangle$ with $n=\dim(V)$, and $a_1,\ldots,a_n\in K^*$. By \cite[Theorem 3.4.1]{joshi-anabelomorphy} one obtains an isomorphism $\alpha:K^*\isom L^*$ and hence $K^{*2}\isom L^{*2}$. Let $(V',Q')$ be the unique quadratic space over $L$ with invariants 
\be
\begin{aligned}
	\dim(V')&=\dim(V),\\
	Q'&=\langle \alpha(a_1),\ldots,\alpha(a_n)\rangle\\
	 \disc(Q')&=\alpha(\disc(Q))=\prod_{i=1}^n \alpha(a_i)\in L^*/L^{*2},\\
	  \varepsilon(V',Q')&=\alpha(\varepsilon(V,Q))=\prod_{i<j}(\alpha(a_i),\alpha(a_j))_L\in\{1,-1\}.
\end{aligned}
\ee
The required correspondence is $(V,Q)/K\mapsto (V',Q')/L$ and inverse correspondence is given by using $\alpha^{-1}:L^*\isom K^*$ by similar formulae.
\ep

\subsection{Aside: Synchronization of Hermitian spaces}
A Hermitian space over $K$ is a triple $(K_1/K,V,\tau)$ where $K_1/K$ is a quadratic extension, $V$ is a $K_1$-vector space and $\tau:V\times V\to K_1$ is a Hermitian pairing with respect to unique non-trivial order two element $\sigma\in\gal(K_1/K)$.
\bpro 
Any anabelomorphism $K\anabmapright{\alpha}L$ provides a synchronization of Hermitian spaces over $K$ and $L$ respectively.
\epro
\bp 
The isomorphism class of $(V,\tau)$ is determined by the quadratic extension $K_1/K$ and two invariants: $\dim(V)$ and $\disc(\tau)\in K^*/N(K_1^*)$ \cite[Chapter 10, Remark 1.4]{scharlau2012-quadratic-book}.
Let $K\anabmapright{\alpha} L$ be an anabelomorphism of $p$-adic fields. Let $H\subset\gamk$ be the open subgroup of index two corresponding to $K_1/K$. Then $H'=\alpha(H)\subset \gaml$ is an open subgroup of index two corresponding to some quadratic extension $L_1/L$. Thus 	$K_1 \anabmapright{\alpha} L_1$. Compatibility with reciprocity maps \cite[Proposition 3.11 and Proposition 4.2]{hoshi-mono} and the anabelomorphism $K_1\anab L_1$, forces the isomorphism of the norm subgroups $N(K_1^*)\isom N(L_1^*)$ of $K^*$ and $L^*$ compatibly with the isomorphism of these groups induced by $\alpha$.  Thus, there is a unique Hermitian space $(L_1/L,V',\tau')$ which is Hermitian with respect to the unique non-trivial element $\sigma'\in\gal(L_1/L)\isom \gal(K_1/K)$ with the following invariants:
\be 
\begin{aligned}
	\dim_L(V')&=\dim_K V,\\
	\disc(\tau')&=\alpha(\disc(\tau))\in L^*/N(L_1^*).
\end{aligned}
\ee
Thus, the Hermitian triple $(K_1/K,V,\tau)$ over $K$ gives rise to the Hermitian triple $(L_1/L,V',\tau')$ with $K_1\anabmapright{\alpha} L_1$. The required correspondence in the reverse direction is given by similar formulae using $\alpha^{-1}$.
\ep

\section{Proof of the synchronization conjecture for a split torus}
\nws
The main conjecture of  \cite{fargues-scholze} is known in the case when $G$ is a torus \cite{zouKonrad2024}. Since \Cref{con:synchronization-conj} is made for split groups, in this section one works with split torii. The explicit description of $\bungk$ given there allows one to prove the existence of $\alpha_{Bun}$ given by \Cref{con:synchronization-conj}.

\bthm\label{thm:anab2}  
Let $L\anabmapright{\alpha} K$ be an anabelomorphism of $p$-adic fields. Then \Cref{con:synchronization-conj} holds for any split torus $T$.
\ethm
\bp 
The required equivalence
\be 
\begin{tikzcd}
	\dbuntl  \ar[r,"\alpha_{Bun}","\isom"'] &  \dbuntk.
\end{tikzcd}
\ee
will be constructed below by constructing an isomorphism of stacks
\be 
\begin{tikzcd}
	\buntl  \ar[r,"\isom"] &  \buntk.
\end{tikzcd}
\ee

Let $G=T=\Gm^r$ (for some integer $r\geq 1$) be a split torus. 
By \cite[Remark 2.1.10]{zouKonrad2024}, one has the following description of the stack
\be \buntk = \coprod_{b\in X_*(T)}[*/\underline{T(K)}],\ee 
where $\underline{T(K)}$ is the condensed abelian group associated with the topological group ${T(K)}$. 

Suppose $L\anabmapright{\alpha} K$ are anabelomorphic $p$-adic fields. Then by \cite[Theorem 3.4.1]{joshi-anabelomorphy}, one has an isomorphism of topological groups
\be\label{eq:KL-mult} 
T(L)=(L^*)^r \isom (K^*)^r=T(K),
\ee
Both these groups are  first countable topological groups, and hence by \cite[Remark 1.6]{scholze2026-condensed} both these groups are also compactly generated topological  abelian groups. From \cite[Proposition 1.7]{scholze2026-condensed} the functor from compactly generated abelian topological groups to condensed abelian groups is fully-faithful, and hence from \eqref{eq:KL-mult}, one obtains a natural isomorphism of condensed groups
\be
\begin{tikzcd}
\underline{T(L)}=\underline{(L^*)}^r\ar [r,"\alpha","\isom"']& \underline{(K^*)}^r=\underline{T(K)}
\end{tikzcd}
\ee
and conversely every  isomorphism between these condensed groups arises from some topological isomorphism \eqref{eq:KL-mult}.
Thus, the anabelomorphism $\alpha$ induces  a natural isomorphism of stacks
\be 
\begin{tikzcd}
	\buntl = \displaystyle{\coprod_{b\in X_*(T)}}[*/\underline{T(L)}] \ar[r,"\alpha_{Bun}","\isom"'] &  \displaystyle{\coprod_{b\in X_*(T)}}[*/\underline{T(K)}] = \buntk.
\end{tikzcd}
\ee
This isomorphism between the stacks $\buntl,\buntk$ gives the equivalence
\be 
\begin{tikzcd}
	\dbuntl  \ar[r,"\alpha_{Bun}","\isom"'] &  \dbuntk.
\end{tikzcd}
\ee

Hence, this proves the existence of $\alpha_{Bun}$ and the theorem. \Cref{thm:anab1} for $G=T$ gives the construction of $\alpha_{Par}$. 

To make the compatibility with $FS_{\psi}$  completely clear, I claim that one has natural isomorphism of structure sheaves
\be\label{eq:z1-isom-torus} 
\O([Z^1(W_K,\td)/\td])\isom \O([Z^1(W_L,\td)/\td]).
\ee
To this one can use \cite[Corollary 4.3.3]{zouKonrad2024}. This asserts that there is a natural isomorphism
\be  
\O([Z^1(W_K,\td)/\td])\isom \invlim_{U\subset T(K)}\Z_\ell[T(K)/U],
\ee
where the inverse limit is over open pro-$p$-subgroups $U\subseteq T(K)$. 

The anabelomorphism $L\anabmapright{\alpha}K$ induces  an homeomorphism \eqref{eq:KL-mult} of topological groups, and being a homeomorphism, it is also an open mapping.  If $U$ runs over open pro $p$-subgroups of $T(L)$, $\alpha(U)$ runs through the open pro $p$-subgroups of $T(K)$. Hence, one has a natural isomorphism
\be 
\Z_\ell[T(L)/U] \mapright{\isom} \Z_\ell[T(K)/\alpha(U)]
\ee
for any open pro $p$-subgroup $U\subset T(L)$ and hence an isomorphism
\be 
\invlim_{U\subset T(L)}\Z_\ell[T(L)/U] \mapright{\isom} \invlim_{U'\subset T(K)}\Z_\ell[T(K)/U']
\ee
in which $U'=\alpha(U)$.
This gives the claimed isomorphism \eqref{eq:z1-isom-torus}.

Now I claim that the anabelomorphism $\alpha$ gives rise to  isomorphic Whittaker datum for $T(K)$ and $T(L)$. This is immediate from \Cref{pr:whittaker}. As $G=T$ is a split torus, its unipotent radical $U=1$ and hence one can also see the correspondence between Whittaker datum directly by establishing the isomorphism of the compact-inductions for $T(L)$ and $T(K)$:
\be  
\begin{tikzcd} 
	\cInd{*}{T(L)} \ar [r,"\alpha","\isom"']& \cInd{*}{T(K)}.
\end{tikzcd}
\ee
But this is immediate from the isomorphism given by \eqref{eq:KL-mult}. In fact, by \cite[Remark 4.2.4]{zouKonrad2024}, 
\be 
\cInd{*}{T(L)} \isom \invlim_{U\subset T(L)}\Z_\ell[T(L)/U]\isom \O([Z^1(W_L,\td)/\td]).
\ee
The correspondence $FS_{\psi}$ (of \cite{fargues-scholze} and for $G=T$ of \cite{zouKonrad2024}) maps the compact induction $\cInd{*}{T(L)}$ i.e. the Whittaker sheaf to the structure sheaf $\O([Z^1(W_L,\td)/\td])$. Thus, this establishes the asserted commutativity of the diagram in \Cref{con:synchronization-conj} and in this case one has 
\be 
\alpha_{Bun}(\sW_L(U_L,\psi'))=\sW_K(U_K,\psi).
\ee
This completes the proof of \Cref{thm:anab2}.
\ep

\bibliography{../../master/masterofallbibs.bib}
\end{document}